\documentclass[a4paper, 12pt]{article}

\usepackage[koi8-r,cp1251]{inputenc}
\usepackage{amsfonts,amssymb,amsmath, hyperref}
\usepackage[final]{epsfig}
\usepackage{graphicx}

\begin{document}

\begin{center}
\textbf{ HAUSDORFF OPERATORS ON  COMPACT ABELIAN GROUPS}
\end{center}

\begin{center}
A. R. Mirotin
\end{center}

\begin{center}
amirotin@yandex.ru
\end{center}

\

{\small Necessary and sufficient conditions are given for boundedness of  Hausdorff  operators on generalized Hardy spaces $H^p_E(G)$, real Hardy space $H^1_{\mathbb{R}}(G)$, $BMO(G)$, and $BMOA(G)$ for compact  Abelian group $G$. Surprisingly, these conditions turned out to be the same for all groups and spaces under consideration. Applications to Dirichlet series are given. The case of the space of continuous functions on $G$ and examples  are also considered.}

\


\

Key wards: Hausdorff  operator, compact  Abelian group, Hardy space, BMO, Dirichlet series.

\

\section{Introduction}

Hausdorff operators are closely connected with classical harmonic analysis (see, e.g., \cite{CFW}, \cite{Ls}, \cite[Chapter XI]{H}, \cite[Section 3]{Boyd}, or \cite{STZ}). The modern stage in the development of this theory  begins with the work by E.~ Liflyand and F.~ M\`{o}ricz \cite{LM}.
The concept of a Hausdorff operator in the general framework of topological groups was introduced by the author in \cite{JMAA} as a generalization of the classical definition in Euclidean spaces \cite{LL}, \cite{BM} and the definition in $p$-adic spaces \cite{Vol} (see Definition 3  below).

In  \cite{JMAA} sufficient conditions  were  given for boundedness of a Hausdorff  operator on the atomic (real)  Hardy space over a  locally compact metrizable  group  that satisfies the so-called doubling property. Generalizations to  homogeneous spaces of Lie groups appeared in \cite{JLieTh}. The case of locally compact groups with local doubling property and their homogeneous spaces  was  considered in \cite{AnMath}.
But there are compact connected Abelian groups that are not metrizable (e.~g., the Bohr compactum $\textbf{b}\mathbb{R}$), or metrizable but without local doubling property (e.~g., the infinite dimensional torus $\mathbb{T}^\infty$\footnote{The author is indebted to Professor A. Bendikov  for this observation.}). The aim of this work is to give necessary and sufficient conditions for boundedness of  Hausdorff  operators on Hardy spaces and $BMO$ for this case, as well. Surprisingly, these conditions   turned out to be the same for all groups and spaces under consideration. The case of the space of continuous functions and examples that may be of interest in their own way are also considered. It should be noted that Hausdorff  operators on Hardy spaces $H^1$ and $BMO$ over Euclidean spaces $\mathbb{R}^d$ were studied in \cite{LL}.

Due to an ingenious identification of
Bohr, a lot of theory of ordinary Dirichlet series may be seen as a sub-theory of Fourier analysis on the infinite dimensional torus. This observation and especially the seminal result of  Hedenmalm,  Lindqvist, and Seip \cite{HLS} give us an opportunity for an
application of our results to ordinary  Dirichlet series. In  the case of general Dirichlet series we use similar results obtained by Defant and Schoolman  in \cite{DS2019}, \cite{I. Schoolmann}, and \cite{DSservey}. Based on this results,  classes of bounded Hausdorff  operators  that act in some  classical spaces of  Dirichlet series (ordinary and general)  are introduced.

\

\section{Preliminaries}

\

This section collects all preliminary information we need in the next parts of the paper.

  In the following unless otherwise stated $G$ stands for a compact and connected Abelian group  with  normalized Haar measure $\nu$ and a total order (which agrees with the group structure) is fixed on its dual group $X$.  In turn, $X$ is the dual group for $G$ by the Pontryagin - van Kampen theorem.
  Let $X_+:=\{\chi\in X:\chi\geq \mathbf{1}\}$  be the positive cone in   $X$ ($\mathbf{1}$ denotes the unit character). In other words,  $X_+$ is a subsemigroup of $X$ such that $X_+^{-1}\cup X_+=X$, and $X_+^{-1}\cap X_+=\{\mathbf{1}\}$  (see, e.g., \cite[Chapter 8]{Rud}). We put also $X_-:=X\setminus X_+$. Then $X_-= X_+^{-1}\setminus \{\mathbf{1}\}$.

    As is well known, a (discrete) Abelian group $X$ can be totally ordered if and
only if it is torsion-free, which in turn is equivalent to the
condition that its character group $G$ is connected. In general the group $X$
may possess many different total orderings.

 In applications, often $X$ is a dense
 subgroup of $\mathbb{R}^d$ endowed with the discrete topology and $G=\mathbf{b}X$ is its Bohr compactification, or $X = \mathbb{Z}^d$ so that
$G = \mathbb{T}^d$ is the $d$-torus ($\mathbb{T}$ is the circle group and $\mathbb{Z}$ is the group of integers). Other interesting examples are the infinite dimensional torus $\mathbb{T}^\infty$ (see Section 6 and Examples 1 and 5 below), an  $\mathbf{a}$-adic solenoids $\Sigma_\mathbf{a}$ (see Example 6 below), and their finite and countable products (see, e.g., \cite{DS2019}).

We denote by $\mathrm{Aut}(H)$  the group of topological automorphisms of a topological group $H$ endowed with its natural topology (see, e.g., \cite{HiR}). If $H$ is Abelian and  $A\in \mathrm{Aut}(H)$   the \textit{dual automorphism} $A^*\in \mathrm{Aut}(\widehat{H})$ is defined by the rule
$$
A^*(\xi):=\xi\circ A,\quad  \xi\in \widehat{H}.
$$
  \cite[(24.37), (24.41)]{HiR}.

  In the following $\mathrm{Aut}_+(X)$ stands for the subset of $\mathrm{Aut}(X)$ that consists of all \textit{ordered automorphisms} of the group $X$ with respect to the given order. By definition, these automorphisms preserve the order (equivalently, these automorphisms map the positive cone $X_+$ into itself).

  We denote also by $\mathrm{Aut}(G)^+$ the set of such $A\in \mathrm{Aut}(G)$  that $A^*\in \mathrm{Aut}_+(X)$.

The following simple lemma will be useful.

\textbf{Lemma 1.} 1) \textit{Let $G$ be locally compact Abelian group.  If $A\in\mathrm{ Aut}(G)$ then $(A^*)^{-1}=(A^{-1})^*$.}

2) \textit{Let $G$ be  compact and connected  Abelian group. Then $\mathrm{ Aut}_+(X)$ is a subgroup of  $\mathrm{ Aut}(X)$.
}

3) \textit{Let $G$ be  compact and connected  Abelian group. Then $\mathrm{ Aut}(G)^+$  is a subgroup of  } $\mathrm{ Aut}(G)$.

Proof. 1) The  map $A\mapsto A^*$ is a topological anti-isomorphism of $\mathrm{Aut}(G)$ onto $\mathrm{Aut}(X)$ \cite[Theorem (26.9)]{HiR}. It follows that $(A^{-1})^*=(A^*)^{-1}$ for   $A\in \mathrm{Aut}(G)$.

 2) We shall show that if $\tau\in\mathrm{ Aut}_+(X)$ then $\tau^{-1}\in\mathrm{ Aut}_+(X)$, as well. Since  $X_-= X_+^{-1}\setminus \{\mathbf{1}\}$, the map $\chi\mapsto \chi^{-1}$ is a bijection of $X_+\setminus \{\mathbf{1}\}$ onto $X_-$. Let us assume that $\tau\in\mathrm{ Aut}_+(X)$ and $\chi\in X_+\setminus \{\mathbf{1}\}$, but $\xi:=\tau^{-1}(\chi)\in X_-$. Then $\xi^{-1}\in X_+\setminus \{\mathbf{1}\}$, and  $\chi=\tau(\xi)=(\tau(\xi^{-1}))^{-1}\in (X_+\setminus \{\mathbf{1}\})^{-1}\subset X_-$, a contradiction.

3) This is an immediate consequence of 1) and 2).

We denote by  $\widehat{\varphi}$
 the Fourier transform of $\varphi\in L^1(G)$, and by $\|\cdot\|_\infty$ the norm in $L^\infty(G)$. We put also
  $$
  \|f\|_p=\left(\int\limits_G|f|^pd\nu\right)^{1/p}
  $$
  for $f\in L^p(G)\ (0< p<\infty)$.

\

In the following the compliment $X\setminus E$ of the subset $E\subset X$ will be denoted by $E^c$.

The next class of spaces is important in particular for general  Hilbert transform \cite{ijpam} and  for the theory of  Dirichlet series (see, e.g., \cite{DS2019} and section 7 below).

\textbf{Definition 1}. \cite{ijpam},  \cite{DS2019}. Let $G$ be compact Abelian group, $1\leq p\leq\infty$, and $E\subset X$
a non-empty set. The \textit{generalized  Hardy space} $H^{p}_E(G)$ is the closed
subspace of  $L^p(G)$ defined as follows
$$
H^{p}_E(G)=\{f\in L^p(G):\widehat{f}(\chi)=0\ \forall\chi\in E^c\}.
$$

The case where $G$ is connected and $E=X_+$ is due to Helson and Lowdenslager (see, e.g., \cite{Rud}). We shall write
$H^p(G)$ instead of $H^{p}_{X_+}(G)$ in this case.
 In particular, $H^2(G)$ is the subspace of $L^2(G)$ with Hilbert basis $X_+$.
We denote by $P_+$  the orthogonal projection $L^2(G)\to H^2(G)$ , and $P_-=I-P_+$.

Of course, the space  $H^p(G)$ (as well as the spaces  $H^p_{\mathbb{R}}(G)$,  $BMO(G)$, and $BMOA(G)$ considered below) depends on the chosen order in $X$.

 For  every $u\in L^2(G,\mathbb{R})$
 there is a unique  $\widetilde{u}\in L^2(G,\mathbb{R})$ such that $\widehat{\widetilde{u}}({\bf 1})=0$ and
$u+{\it i}\widetilde{u}\in H^2(G)$.
The linear continuation of the mapping
$u\mapsto \widetilde{u}$ to the complex $L^2(G)$ is called a {\it Hilbert transform} on
$G$. This operator extends to a  bounded operator $\varphi\mapsto \widetilde\varphi$ on $L^p(G)$ for $1<p<\infty$ (generalized Marcel Riesz's inequality), in particular  $\|\widetilde{\varphi}\|_2\leq\|\varphi\|_2$ for every $\varphi\in L^2(G)$  \cite[8.7]{Rud}, \cite[Theorem 8, Corollary 20]{ijpam}. Let $\mathcal{F}f=\widehat{f}$ be the Fourier transform on $G$. Then the next formula holds
$$
\widehat{\widetilde{f}}=-i\mathrm{sgn}_{X_+}\widehat{f},
$$
where $\mathrm{sgn}_{X_+}(\chi)=1$ for $\chi\in X_+\setminus\{\mathbf{1}\}$, $\mathrm{sgn}_{X_+}(\mathbf{1})=0$, and $\mathrm{sgn}_{X_+}(\chi)=-1$ for $\chi\in X\setminus X_+$ \cite{ijpam}.

 Note also that  the Hilbert transform is a continuous map from $L^1(G)$ to $L^p(G)$ for $0<p<1$ (see, e.~g., \cite[Theorem 8.7.6]{Rud}).

\textbf{Definition 2} \cite{Trudy}. We define the space $BMO(G)$ of
\textit{functions of bounded mean oscillation on   $G$} and its subspace $BMOA(G)$,
 as follows
$$
BMO(G):=\{f+\widetilde{g}: f,g\in L^{\infty}(G)\},
$$
$$
BMOA(G):=BMO(G)\cap H^1(G),
$$
$$
\|\varphi\|_{BMO}:=\inf\{\|f\|_{\infty}+\|g\|_{\infty}:\varphi=f+\widetilde{g},
f,g\in L^{\infty}(G)\}
$$
for $\varphi\in BMO(G)$.

 \textbf{Lemma 2.} \cite[Lemma 1]{Indag}.
\textit{The following equalities hold:}

$
1)\  BMO(G)=P_-L^{\infty}(G)+P_+L^{\infty}(G),
$
\textit{with an equivalent norm}
$$
\|\varphi\|_{\ast}:=\inf\{\max(\|f_1\|_\infty, \|g_1\|_\infty):\varphi=P_-f_1+P_+g_1, f_1, g_1\in L^\infty(G)\};
$$

2) \ $BMOA(G)=P_+L^{\infty}(G)$. \textit{Moreover, for the norm}
$$\|\varphi\|_{\ast}=\inf\{\|h\|_\infty:\varphi=P_+h,\
h\in L^\infty(G)\}
$$
\textit{in this space the following inequalities take place: \footnote{Here we correct a typo made in \cite[p. 139]{Trudy}.}}
$$\frac{2}{3}\|\varphi\|_
{BMO}\leq \|\varphi\|_
{\ast}\leq 2\|\varphi\|_
{BMO}.$$

\textbf{Definition 3}. \cite{Indag}  We define the space $H^1_{\mathbb{R}}(G)$ (the \textit{real $H^1$ space on $G$})
as the completion of the space ${\rm Pol}(G,\mathbb{R})$ of real-valued trigonometric polynomials on $G$
 with respect to the norm
 $$
 \|q\|_{1\ast}:=\|P_-q\|_1+\|P_+q\|_1.
 $$
 We denote the norm in  $H^1_{\mathbb{R}}(G)$ by  $\|\cdot\|_{1\ast}$, too.

 The notation $H^1_{\mathbb{R}}(G)$ should not lead to the confusion with $H^p_E(G)$ from the Definition 1.

\textbf{Lemma 3} \cite[Proposition 1]{Indag}. (i) \textit{Projectors $P_{\pm}$, and the Hilbert transform  are  bounded operators on } $H^1_{\mathbb{R}}(G)$;

 (ii) \textit{restrictions $P_{\pm}|{\rm Pol}(G,\mathbb{R})$ extend to  bounded operators $P_{\pm}^1$ from $ H^1_{\mathbb{R}}(G)$ to $L^1(G)$ and
 $$
 \|f\|_{1*}=\|P_-f\|_{1*}+\|P_+f\|_{1*}=\|P_-^1f\|_{1}+\|P_+^1f\|_{1}\    (f\in H^1_{\mathbb{R}}(G));
 $$
 }

 (iii) $ H^1_{\mathbb{R}}(G)={\rm Im}P_-\dotplus {\rm Im}P_+$ (\textit{the direct sum of closed subspaces});

 (iv) $\cup_{p>1}L^p(G,\mathbb{R})\subset  H^1_{\mathbb{R}}(G) \subset L^1(G,\mathbb{R})$;

 (v) $\|f\|^{\sim}:=\|f\|_1+\|\widetilde{f}\|_1$ \textit{is an equivalent norm in} $ H^1_{\mathbb{R}}(G)$;

 (vi)  $H^1_{\mathbb{R}}(G)={\rm Re}H^1(G)$.

In \cite{JMAA} the next definition was proposed.

 \textbf{Definition 4} \cite{JMAA}. Let  $(\Omega,\mu)$ be  a measure space, $G$ a topological group, $A:\Omega\to \mathrm{Aut}(G)$ a measurable map,
   and $\Phi$ a locally $\mu$-integrable function on $\Omega$.
   We define the \textit{Hausdorff  operator} with the kernel $\Phi$ over the group $G$  by the formula
$$
(\mathcal{H}_{\Phi,A}f)(x)=\int_\Omega \Phi(u)f(A(u)(x))d\mu(u).
$$

In particular, we get a class of \textit{discrete Hausdorff  operators} of the form
$$
f\mapsto\sum_{u\in \Omega} \Phi(u)(f\circ A(u))
$$
where $\Omega$ is a countable set endowed with the counting measure.

Throughout  we denote by $\mathcal{L}(Y)$ the space of linear bounded operators on a
 normed space $Y$.

By \cite[Lemma 1]{JMAA}  an operator $\mathcal{H}_{\Phi,A}$ is bounded on $L^p(G)$ ($1\le p\le\infty$) for a locally compact group  $G$ provided $\Phi(u)(\mathrm{mod} A(u))^{-1/p}\in L^1(\Omega,\mu)$, and
$$
\|\mathcal{H}_{\Phi,A}\|_{\mathcal{L}(L^p(G))}\le \int_\Omega |\Phi(u)|(\mathrm{mod} A(u))^{-1/p}d\mu(u). \eqno(1)
$$

\textbf{Example 1.} Let $\mathbb{T}^\infty$ be the infinite-dimensional torus (the product of a countably many copies of the circle group) and $\mathcal{C}:=\{-1,1\}^\infty$   a Cantor group endowed by some regular Borel measure $\mu$
(e.g., $\mu$ is the normalized Haar measure of the compact group  $\mathcal{C}$).  The  group $\mathcal{C}$
 acts on $\mathbb{T}^\infty$ by coordinate-wise  automorphisms  $A(u)(x)=x^u:=(x_j^{u_j})_{j\in\mathbb{ N}}$ where $u=(u)_{j\in\mathbb{ N}}$, $u_j\in \{-1,1\}$, and $x=(x_j)_{j\in\mathbb{ N}}$, $x_j\in \mathbb{T}$. Thus, we get a Hausdorff operator
 $$
 \mathcal{H}_{\Phi}f(x)=\int_{\mathcal{C}}\Phi(u)f(x^u)d\mu(u).
 $$
Since  $\mathbb{T}^\infty$  is unimodular, $\mathrm{mod} A(u)=1$ and so this operator is bounded on $L^p(\mathbb{T}^\infty)$ ($1\le p\le\infty$) for $\Phi\in L^1(\mu)$ and $\|\mathcal{H}_\Phi\|_{\mathcal{L}(L^p(\mathbb{T}^\infty))}\le \|\Phi\|_{L^1(\mu)}.$

\

\section{Commuting Relations for Hausdorff  Operator}

In this section we shall show that Hausdorff  operator commutes in some sense both with the Fourier transform and the Hilbert transform.

%Despite being mainly interested in the case of compact groups, we write this part for general LCA groups.

\textbf{Theorem 1} (cf. \cite[Theorem 4.4]{Ls}).
 (i) \textit{Let  $G$ be compact (not necessary  connected) Abelian  group,
 $f\in L^1(G)$, and  $\Phi\in L^1(\mu)$. Then}
$$
(\mathcal{H}_{\Phi,A}f)^{\wedge}=\mathcal{H}_{\Phi,(A^*)^{-1}} \widehat{f}.
$$

(ii)\textit{ Let  $G$ be compact and  connected Abelian group, $f\in L^2(G)$,  $\Phi\in L^1(\mu)$,  and    $A(u)\in \mathrm{Aut}(G)^+$  for $\mu$-a.~e. $u\in \Omega$. Then}
$$
\mathcal{H}_{\Phi,A}\widetilde{f}=(\mathcal{H}_{\Phi,A}f)^{\thicksim}.
$$

Proof. (i) By the Fubini theorem
$$
(\mathcal{H}_{\Phi,A}f)^{\wedge}(\chi)=\int_G\left(\int_\Omega \Phi(u)f(A(u)(x)d\mu(u)\right)\overline{\chi(x)}d\nu(x)
$$
$$
=\int_\Omega \Phi(u)\left(\int_G f(A(u)(x)\overline{\chi(x)}d\nu(x)\right)d\mu(u).
$$
 Moreover, since $G$ is unimodular, we have $\mathrm{mod} A(u)=1$, and we get putting $y=A(u)(x)$ that
$$
\int_G f(A(u)(x)\overline{\chi(x)}d\nu(x)=\int_G f(y)\overline{\chi(A(u)^{-1}(y))}d\nu(y)=\widehat{f}((A(u)^*)^{-1}(\chi)).
$$
So,
$$
(\mathcal{H}_{\Phi,A}f)^{\wedge}=\mathcal{H}_{\Phi,(A^*)^{-1}}(f^{\wedge}).
$$

(ii) Note that $\widetilde{f}\in L^2(G)$. Then in view of (i) one has for all $\chi\in X$ that
$$
\mathcal{F}(\mathcal{H}_{\Phi,A}\widetilde{f})(\chi)=\mathcal{H}_{\Phi,(A^*)^{-1}} \widehat{\widetilde{f}}(\chi)=(\mathcal{H}_{\Phi,(A^*)^{-1}}(-i\mathrm{sgn}_{X_+} \widehat{f}))(\chi)
$$
$$
=-i\int_\Omega \Phi(u)\mathrm{sgn}_{X_+}((A^*)^{-1}(\chi))\widehat{f}((A(u)^*)^{-1}(\chi))d\mu(u).
$$
Since $(A(u)^*)^{-1}$ is an order automorphism for $\mu$-a.~e. $u\in \Omega$, one has $\mathrm{sgn}_{X_+}((A(u)^*)^{-1}(\chi))=\mathrm{sgn}_{X_+}(\chi)$ a.~e. This yields  (again by (i)) that
$$
\mathcal{F}(\mathcal{H}_{\Phi,A}\widetilde{f})(\chi)=-i \mathrm{sgn}_{X_+}(\chi)\int_\Omega \Phi(u)\widehat{f}((A(u)^*)^{-1}(\chi))d\mu(u)
$$
$$
=-i\mathrm{sgn}_{X_+}(\chi)\mathcal{F}(\mathcal{H}_{\Phi,A}f)(\chi)=\mathcal{F} (\mathcal{H}_{\Phi,A}f)^{\sim}(\chi),
$$
which completes the proof.

\textbf{Corollary 1.} \textit{Let $\Phi\in L^1(\mu)$ and $A(u)\in \mathrm{Aut}(G)^+$  for $\mu$-a.~e. $u\in \Omega$. Then the range of  $\mathcal{H}_{\Phi,A}$ in the space $L^2(G)$ is invariant with respect to the Hilbert transform.}

\

\section{ Hausdorff Operators on  Spaces $H^{p}_E(G)$ and $BMOA(G)$}

The next theorem deals with general compact Abelian groups.

\textbf{Theorem 2.}  \textit{Let $G$ be compact  (not necessary connected) Abelian group, $E\subset X$, and $(A(u)^*)^{-1}:E^c\to E^c$ for $\mu$-a.~e. $u\in \Omega$.
 The Hausdorff operator $\mathcal{H}_{\Phi,A}$ is bounded on $H^{p}_{E}(G)$ ($1\le p\le\infty$) if $\Phi\in L^1(\mu)$. In this case,}
$$
\|\mathcal{H}_{\Phi,A}\|_{\mathcal{L}(H^{p}_{E})}\le \|\Phi\|_{L^1(\mu)}.
$$

Proof.  Let $\Phi\in L^1(\mu)$. Since $G$ is unimodular, $\mathrm{mod} A(u)=1$. Thus, as was mentioned in the Introduction, the operator $\mathcal{H}_{\Phi,A}$ is bounded in $L^p(G)$ and formula  (1) holds with  $\mathrm{mod} A(u)=1$. So, it suffices to show that  $\mathcal{H}_{\Phi,A}$ acts in $H^{p}_{E}(G)$. In other wards, it suffices to show that for each $f\in H^{p}_{E}(G)$ the Fourier transform
of $\mathcal{H}_{\Phi,A}f$ is concentrated on $E$.   But by the Theorem 1 (i)
$$
(\mathcal{H}_{\Phi,A}f)^{\wedge}(\chi)=  \int_\Omega \Phi(u)\widehat{f}((A(u)^*)^{-1}(\chi))d\mu(u).
$$
    Let $\chi\in E^c$.  Since $\widehat{f}$ is concentrated on $E$, we have $\widehat{f}((A(u)^*)^{-1}(\chi))=0$ for $\mu$-a.~e. $u\in \Omega$.  It follows that $(\mathcal{H}_{\Phi,(A)}f)^{\wedge}(\chi)=0$, too.
This completes the proof.

\textbf{Corollary 2. } \textit{Let $G$ be compact and connected Abelian group, and $A(u)\in \mathrm{Aut}(G)^+$ for $\mu$-a.~e. $u\in \Omega$.
 The Hausdorff operator $\mathcal{H}_{\Phi,A}$ is bounded on $H^p(G)$ ($1\le p\le\infty$) if and only if $\Phi\in L^1(\mu)$. In this case,}
$$
\|\mathcal{H}_{\Phi,A}\|_{\mathcal{L}(H^p)}\le \|\Phi\|_{L^1(\mu)}.
$$

Proof. Since $\mathbf{1}\in H^p(G)$, the "only if" part is obvious. Now let $\Phi\in L^1(\mu)$.  In our case $E=X_+$. So, it suffices to show that $(A(u)^*)^{-1}:X_-\to X_-$ for $\mu$-a.~e. $u\in \Omega$.  Indeed, let $\chi\in X_ -=X\setminus X_+$. Then $\chi^{-1}\in X_+\setminus\{\mathbf{1}\}$
and therefore $(A(u)^*(\chi))^{-1}=A(u)^*(\chi^{-1})\in X_+\setminus\{\mathbf{1}\}$. Thus,  $A(u)^*(\chi)\in X\setminus X_+$.
This completes the proof.

From now on, we denote  by $Y^*$ the dual of the space $Y$ and by $B^*$ the adjoint
of an operator $B\in \mathcal{L}(Y)$.

 For the proof of our next theorem we need the following

\textbf{ Theorem A}. (\cite{Indag}, Theorem 1). \textit{
For every $\varphi\in BMOA(G)$ the formula
$$
F_\varphi(f)=\int\limits_G f\overline \varphi d\nu
$$
defines a linear functional on  $H^\infty(G)$,
and this functional extends uniquely to a continuous linear functional $F_\varphi$  on $H^1(G)$. Moreover, the correspondence  $\varphi\mapsto F_\varphi$ is an isometrical isomorphism of $(BMOA(G),\|\cdot\|_*)$  and $H^1(G)^*$, and a topological isomorphism of $(BMOA(G),\|\cdot\|_{BMO})$  and $H^1(G)^*$.}

\textbf{Theorem 3.} \textit{Let  $A(u)\in \mathrm{Aut}(G)^+$ for $\mu$-a.~e. $u\in \Omega$.
 The Hausdorff operator $\mathcal{H}_{\Phi,A}$  is bounded on the  space $BMOA(G)$ if and only if $\Phi\in L^1(\mu)$.
Moreover,}
$$
\|\mathcal{H}_{\Phi,A}\|_{\mathcal{L}(BMOA)}\le\|\Phi\|_{L^1(\mu)}.
$$

Proof. Since $\mathbf{1}\in BMOA(G)$, the "only if" part is obvious. Now  let  $\Phi\in L^1(\mu)$.
In view of Theorem 2 for the proof it suffices to show that $\mathcal{H}_{\Phi,A}=\mathcal{H}_{\overline{\Phi},A^{-1}}^*$ where $\mathcal{H}_{\overline{\Phi},A^{-1}}$ is considered in $ H^1(G)$. To this end we shall employ Theorem A. Let $f\in H^\infty(G)$. Then it is clear that  $\mathcal{H}_{\overline{\Phi},A^{-1}}f\in H^\infty(G)$, too and for every $\varphi\in BMOA(G)$ we have
$$
\mathcal{H}_{\overline{\Phi},A^{-1}}^*(F_\varphi)(f):=F_\varphi(\mathcal{H}_{\overline{\Phi},A^{-1}}f)=
\int_G\left(\int_\Omega \overline{\Phi(u)}f(A(u)^{-1}(x)d\mu(u)\right)\overline{\varphi(x)}d\nu(x)
$$
$$
=\int_\Omega \overline{\Phi(u)}\left(\int_G f(A(u)^{-1}(x)\overline{\varphi(x)}d\nu(x)\right)d\mu(u).
$$
 by the Fubini theorem.

Further, as in the proof of Theorem 1, we get putting $y=A(u)(x)$ that
$$
\int_G f(A(u)^{-1}(x)\overline{\varphi(x)}d\nu(x)=\int_G f(y)\overline{\varphi(A(u)(y))}d\nu(y).
$$
Thus, (again by the Fubini theorem)
$$
\mathcal{H}_{\overline{\Phi},A^{-1}}^*(F_\varphi)(f)=\int_G f(y)\left(\int_\Omega \overline{\Phi(u)} \overline{\varphi(A(u)(y))}d\mu(u)\right)d\nu(y)
$$
$$
=\int_G f(y)\overline{ \mathcal{H}_{\Phi,A}\varphi(y)}d\nu(y)=F_\psi(f),
$$
where $\psi=\mathcal{H}_{\Phi,A}\varphi$. Since by Theorem A every  continuous linear functional   on $H^1(G)$ is uniquely defied by its values on $H^\infty$, it follows that  $\mathcal{H}_{\overline{\Phi},A^{-1}}^*(F_\varphi)=F_\psi$. If we identify (again by Theorem A)
$F_\varphi$ with $\varphi$ and $F_\psi$ with $\psi$ we have
$$
\mathcal{H}_{\Phi,A}\varphi=\mathcal{H}_{\overline{\Phi},A^{-1}}^*\varphi,
$$
which completes the proof.

For the next corollary we need the following

\textbf{Definition 5.} \cite{Rud}. We call a subset $E\subset X_+$   \textit{lacunary} (in the sense of Rudin) if there is a constant  $K_E$ such that  the number of terms of  the set $\{\xi\in E:\chi\leq\xi\leq\chi^2\}$ do not exceed $K_E$ for every $\chi\in X_+$.

\textbf{Corollary 3.}  \textit{Let the subset $E$ of $X_+$ be lacunary, $A(u)\in \mathrm{Aut}(G)^+$ for $\mu$-a.~e. $u\in \Omega$,
 and $\Phi\in L^1(\mu)$.
%$(A(u)^*)^{-1}\in \mathrm{Aut}(G)^+$ for $\mu$-a.~e. $u\in \Omega$.
Then $\mathcal{H}_{\Phi,A}$ is a bounded  operator from $H^2_E(G)$ into $(BMOA(G),\|\cdot\|_{BMO})$ and}
$$
\|\mathcal{H}_{\Phi,A}\|_{H^2_E\to BMOA}\le 3\sqrt{K_E}\|\Phi\|_{L^1}.
$$

Proof. Let $Pol_E(G):=\mathrm{span}_{\mathbb{C}}(E)$ be the space of $E$-polynomials.
It is known \cite[Propositiuon 3.14]{DS2019}, \cite[Lemma 1]{MatSb} that $Pol_E(G)$ is a dense  subspace of $H^p_E(G)$ for all $p\in [1,\infty)$.
Since $E\subset X_+$, we have  $H^p_E(G)\subset H^p(G)$.
Let $\varphi\in Pol_E(G)$. Then $\varphi\in H^1_E(G)\cap H^2(G)$ and by \cite[Theorem 3]{Indag} one has  $\varphi\in BMOA(G)$ and $\|\varphi\|_{BMO}\le  3\sqrt{K_E}\|\varphi\|_{H^2}$. Now Theorem 3 yields, that
$$
\|\mathcal{H}_{\Phi,A}\varphi\|_{BMO}\le \|\Phi\|_{L^1}\|\varphi\|_{BMO}\le  3\sqrt{K_E}\|\Phi\|_{L^1}\|\varphi\|_{H^2}
$$
and the result follows.

In conclusion to this section, we discuss the necessity of the condition in Corollary 2 of Theorem 2 and Theorem  3.

\textbf{Proposition 1. }  \textit{Let $G$ be metrizable. Assume  in addition to the assumptions of Definition 4 that $\int_{E}\Phi d\mu\ne 0$ for every measurable $E\subset \Omega$, $\mu(E)>0$. If the Hausdorff operator
$\mathcal{H}_{\Phi,A}$ acts in $H^1(G)$ or $BMOA(G)$ then $A(u)\in \mathrm{Aut}(G)^+$ for a.~e. $u\in \Omega$.}

Proof. Since $X_+\subset BMOA(G)\subset H^1(G)$, we have for every $\chi\in X_+$ and every $\xi\in X_-=X\setminus X_+$  by Theorem 1 that
$$
(\mathcal{H}_{\Phi,A}\chi)^{\wedge}(\xi)=(\mathcal{H}_{\Phi,(A^*)^{-1}} \widehat{\chi})(\xi)=0.
$$
On the other hand, the orthogonality of characters of $G$ implies that $\widehat{\chi}=1_{\{\chi\}}$, where $1_A$ stands for the indicator function of a subset $A$ of $X$. Thus,
$$
0=(\mathcal{H}_{\Phi,(A(u)^*)^{-1}} 1_{\{\chi\}})(\xi)=\int_{E(\chi,\xi)}\Phi(u)d\mu(u),
$$
where
$$
E(\chi,\xi)=\{u\in \Omega: (A(u)^*)^{-1}(\xi)=\chi\}=\{u\in \Omega: A(u)^*(\chi)=\xi\}.
$$
 Therefore $\mu(E(\chi,\xi))=0$ for an   arbitrary $\chi\in X_+$ and  $\xi\in X_-$. Moreover,
 $$
 \{u\in \Omega: A(u)^*:X_+\nrightarrow X_+\}=\cup\{E(\chi,\xi): \chi\in X_+, \xi\in X_-\}.
 $$
  Since $G$ is metrizable, $X$ is countable (see, e.g., \cite[Corollary of Theorem 29]{Morris}).  It follows that $A(u)^*:X_+\to X_+$ for $\mu$-a.~e. $u$, which completes the proof.

\

\section{Hausdorff Operators on  Spaces $BMO(G)$ and  $H^1_\mathbb{R}(G)$}

\textbf{Theorem 4.} \textit{Let  $A(u)\in \mathrm{Aut}(G)^+$ for $\mu$-a.~e. $u\in \Omega$.
 The Hausdorff operator $\mathcal{H}_{\Phi,A}$ is bounded on $BMO(G)$ if
and only if $\Phi\in L^1(\mu)$. In this case,}
$$
\|\mathcal{H}_{\Phi,A}\|_{\mathcal{L}(BMO)}\le \|\Phi\|_{L^1(\mu)}.
$$

Proof. The necessity  is obvious. Let $\Phi\in L^1(\mu)$. Every function $\varphi\in BMO(G)$ has the form  $\varphi=f+\widetilde{g}$ where $f,g\in L^{\infty}(G)$. Then by Theorem 1
$$
\mathcal{H}_{\Phi,A}\varphi=\mathcal{H}_{\Phi,A} f+(\mathcal{H}_{\Phi,A} g)^{\sim}.
$$
Note that $\mathcal{H}_{\Phi,A} f, \mathcal{H}_{\Phi,A} g\in L^{\infty}(G)$. Thus,
$$
\|\mathcal{H}_{\Phi,A}\varphi\|_{BMO}\le \|\mathcal{H}_{\Phi,A} f\|_{\infty}+\|\mathcal{H}_{\Phi,A} g\|_{\infty}\le \|\Phi\|_{L^1(\mu)}(\|f\|_{\infty}+\|g\|_{\infty}),
$$
and the result follows.

Below we shall use the following

\textbf{ Theorem B} (\cite{Indag}, Theorem 2). \textit{
For every  $\varphi\in BMO(G,\mathbb{R})$ the linear functional
$$
F_\varphi(q)=\int\limits_G q \varphi d\nu \eqno(2)
$$
 on  ${\rm Pol}(G,\mathbb{R})$ extends uniquely to a continuous linear functional $F_\varphi$  on $H^1_{\mathbb{R}}(G)$. Moreover, the correspondence  $\varphi\mapsto F_\varphi$ is an isometrical isomorphism of $(BMO(G,\mathbb{R}),\|\cdot\|_*)$  and $H^1_{\mathbb{R}}(G)^*$, and  a topological isomorphism of  $(BMO(G,\mathbb{R})$, $\|\cdot\|_{BMO})$  and $H^1_{\mathbb{R}}(G)^*$.}

\textbf{Corollary 4.}\textit{ Theorem B is valid with $q\in L^2(G,\mathbb{R})$ in place of $q\in \mathrm{Pol}(G,\mathbb{R})$.}

Proof. Since ${\rm Pol}(G,\mathbb{R})\subset L^2(G,\mathbb{R})$ and ${\rm Pol}(G,\mathbb{R})$ is dense in $H^1_{\mathbb{R}}(G)$,  it suffices to show that the right-hand side in (2) is continuous on the set  $L^2(G,\mathbb{R})$ with respect to the $H^1_{\mathbb{R}}(G)$ norm. Let (Lemma 2) $\varphi=P_-g+P_+h$, where $g, h\in L^\infty(G)$.
 Then for every $q\in L^2(G,\mathbb{R})$ one  has that ($q$ is real valued)
$$
\int_G q\varphi d\nu= \int_G P_-gq d\nu+\int_G P_+h q d\nu=\int_G g\overline{P_-q d}\nu+\int\limits_G h\overline{P_+q} d\nu.
$$
This yields that
$$
\left|\int_G q\varphi d\nu\right|\leq\max(\|g\|_\infty,\|h\|_\infty)(\|P_-q\|_1+\|P_+q\|_1).
$$
So, $\left|\int_G q\varphi d\nu \right|\leq \|\varphi\|_*\|q\|_{1\ast}$
(we used Lemma 3 and the fact that $P_{\pm}^1 q=P_{\pm}q$ for $q\in L^2(G,\mathbb{R})$)
 and the proof is complete.

\textbf{Theorem 5.} \textit{Let  $A(u)\in \mathrm{Aut}(G)^+$ for $\mu$-a.~e. $u\in \Omega$.
 The Hausdorff operator $\mathcal{H}_{\Phi,A}$ with real valued $\Phi$ is bounded on the real Hardy space $H^1_\mathbb{R}(G)$ if and only if $\Phi\in L^1(\mu)$. Moreover,}
$$
\|\mathcal{H}_{\Phi,A}\|_{\mathcal{L}(H^1_\mathbb{R})}\le \|\Phi\|_{L^1(\mu)}.
$$

Proof. As above, the "only if" part is obvious. Let $\Phi\in L^1(\mu)$. We shall employ Theorem 2 and the fact that $H^1_\mathbb{R}(G)=\mathrm{Re }H^1(G)$ (Lemma 3). Let $g\in  H^1_\mathbb{R}(G)$. Then $g=f+\overline{f}$ where $f\in H^1(G)$. But since $\Phi$ is real, we have
$$
\overline{\mathcal{H}_{\Phi,A}f(x)}=\int_\Omega \Phi(u)\overline{f(A(u)(x)}d\mu(u).
$$
Since $\mathcal{H}_{\Phi,A}$  is linear in $L^1(G)$, it follows that
$$
\mathcal{H}_{\Phi,A}g=\mathcal{H}_{\Phi,A}f+\mathcal{H}_{\Phi,A}\overline{f}=\mathcal{H}_{\Phi,A}f+\overline{\mathcal{H}_{\Phi,A}f}
\in \mathrm{Re }H^1(G)=H^1_\mathbb{R}(G).
$$
Thus, $\mathcal{H}_{\Phi,A}$ acts in $H^1_\mathbb{R}(G)$.
Now we shall apply the closed graph  theorem.
Let $f_n\to f$ and $\mathcal{H}_{\Phi,A}f_n\to g$ in $H^1(G)$. Since there is a continuous embedding $H^1_\mathbb{R}(G)\subset L^1(G,\mathbb{R})$ (Lemma 3), it follows that $f_n\to f$ and $\mathcal{H}_{\Phi,A}f_n\to \mathcal{H}_{\Phi,A}f$ in $L^1(G)$. Thus, $g= \mathcal{H}_{\Phi,A}f$ and the proof of the continuity of $\mathcal{H}_{\Phi,A}$  is complete.

Finally,  due to Corollary 4 as in the proof of Theorem 3 we have $\mathcal{H}_{\Phi,A}^*=\mathcal{H}_{\Phi,A^{-1}}$,
where $\mathcal{H}_{\Phi,A^{-1}}$ is considered in $BMO(G)$. Then by Theorem 4
$$
\|\mathcal{H}_{\Phi,A}\|_{\mathcal{L}(H^1_\mathbb{R})}=\|\mathcal{H}_{\Phi,A^{-1}}\|_{\mathcal{L}(BMO)}\le \|\Phi\|_{L^1(\mu)}.\eqno(3)
$$

\textbf{Remark 1.} It is clear that $\mathbf{1}\in H^1_\mathbb{R}(G)$  and $\|\mathbf{1}\|_{1\ast}=1$. If $\Phi\ge 0$, we have $\mathcal{H}_{\Phi,A}\mathbf{1}=\|\Phi\|_{L^1(\mu)}\mathbf{1}$. Thus, $\|\mathcal{H}_{\Phi,A}\|_{\mathcal{L}(H^1_\mathbb{R})}=\|\Phi\|_{L^1(\mu)}$. Then formula (3) shows that $\|\mathcal{H}_{\Phi,A}\|_{\mathcal{L}(BMO)}=\|\Phi\|_{L^1(\mu)}$, as well. For  $\Phi\ge 0$ similar equalities hold for the spaces $H^p(G)$ ($1\le p\le\infty$), $H^1_{\mathbb{R}}(G)$, and $BMOA(G)$.

\

\section{On the Action of $\mathcal{H}_{\Phi,A}$ in $C(G)$}

The next simple proposition gives  sufficient conditions for boundedness of a Hausdorff operator in $C(G)$.

\textbf{Proposition 2.} \textit{Let  $G$ be compact (not necessary  connected) Abelian  group, and one of the following two conditions holds:}

1) $G$ \textit{is  first-countable};

2) $\Omega$ \textit{is a completely regular topological space with a bounded Radon measure $\mu$, $\Phi$ is a bounded and continuous function on $\Omega$, and the map $\Omega\times G\to G$, $(u,x)\mapsto A(u)(x)$ is continuous.}

\textit{Then $\mathcal{H}_{\Phi,A}$ acts  in the space  $C(G)$  and is bounded if and only if   $\Phi\in L^1(\mu)$ and in this case  $\|\mathcal{H}_{\Phi,A}\|\le \|\Phi\|_{L^1(\mu)}$.}

Proof. The necessity in obvious.  In the case 1) the sufficiency  follows from the Lebesgue theorem, and in the case 2) this follows, e.~g., from \cite[Chapter IX, \S 5, Corollary of Proposition 13]{Burb}.

The following example shows that the conditions of the  previous Proposition are essential, because \textit{in general  $\mathcal{H}_{\Phi,A}$ does not act  in   $C(G)$.
}

\textbf{Example 2.}  Let $G=\mathbf{b}\mathbb{R}$ be the Bohr compactification of the reals  (see, e.~g., \cite[Section 1.8]{Rud}). This means that $G$ is the dual group of the additive group $X:=\mathbb{R}_d$ where the group $\mathbb{R}$ of reals  is endowed with the discrete topology and the usual order. Then $X$ is the dual group of $\mathbf{b}\mathbb{R}$ by the Pontryagin - van Kampen theorem.  The map $\tau_u(\gamma):=u\gamma$ belongs to $\mathrm{Aut}(X)$ for every $u\in \mathbb{R}, u\ne 0$.

 For each $t\in \mathbb{R}$ let $\widehat{t}(\gamma)=e^{-it\gamma}$ be the corresponding continuous character of  $\mathbb{R}$ ($\gamma\in \mathbb{R}$). Then the map $\beta:\mathbb{R}\to \mathbf{b}\mathbb{R}$,  $t\mapsto \widehat{t}$ is a continuous isomorphism of  $\mathbb{R}$ onto a dense subgroup of  $\mathbf{b}\mathbb{R}$ (see, e.~g., \cite[1.8.2]{Rud}).
 So  we identify $\widehat{t}$ with $t\in \mathbb{R}$ and consider $\mathbb{R}$ as a dense subgroup of  $\mathbf{b}\mathbb{R}$.

The space $AP(\mathbb{R})$ of uniformly almost periodic functions on $\mathbb{R}$ (endowed with the $\mathrm{sup}$ norm) is isometrically isomorphic to  $C(\mathbf{b}\mathbb{R})$ via the restriction map  $C(\mathbf{b}\mathbb{R})\to AP(\mathbb{R})$, $g\mapsto g|\mathbb{R}$ (see, e.~g., \cite[1.8.4]{Rud}, \cite[Chapter VIII, \S 41]{Loomis}).

Let $\Omega= \mathbb{R}$, $d\mu(u)=du$, $\Phi\in L^1(\mathbb{R})$. If we assume that the Hausdorff operator
$$
 \mathcal{H}_{\Phi,\tau_u^*}g(x)=\int_{\mathbb{R}} \Phi(u)g(\tau_u^*(x))du
$$
acts  in $C(\mathbf{b}\mathbb{R})$ then
the  operator
$$
 \mathcal{H}_{\Phi}f(t):=\int_{\mathbb{R}} \Phi(u)f(ut)du
$$
acts  in $AP(\mathbb{R})$.

For the proof it suffices to show that
$\mathcal{H}_{\Phi}f= (\mathcal{H}_{\Phi,\tau_u^*}g)|\mathbb{R}$, where $g\in C(\mathbf{b}\mathbb{R})$, $f:=g|\mathbb{R}$. But for $t\in \mathbb{R}$ one has
$$
\tau_u^*(\widehat{t})(\gamma)=\widehat{t}(u\gamma)=e^{-iut\gamma}=\widehat{ut}(\gamma)\  \ (\gamma\in \mathbb{R}).
$$
Thus, $\tau_u^*(\widehat{t})=\widehat{ut}$.  It follows that for $t\in \mathbb{R}$

$$
 \mathcal{H}_{\Phi,\tau_u^*}g(t)= \mathcal{H}_{\Phi,\tau_u^*}g(\widehat{t})=\int_{\mathbb{R}} \Phi(u)g(\tau_u^*(\widehat{t}))du=\int_{\mathbb{R}} \Phi(u)f(\widehat{ut})du=\int_{\mathbb{R}} \Phi(u)f(ut)du
$$
(recall that  we identify $\widehat{t}$ with $t\in \mathbb{R}$), which completes the proof.

In particular, taking  $f(t)=e^{-it}$, we get that for  $\Phi\in L^1(\mathbb{R})\cap C(\mathbb{R})$ the Fourier transform $\widehat{\Phi}$ belongs to $AP(\mathbb{R})$.
But it is known (see, e.~g., \cite[Theorem 3]{Eberlain}) that in this case the measure $\Phi(u)du$ should be discrete, and we get a contradiction.

\

\section{Applications to Dirichlet Series}

\subsection{Ordinary Dirichlet Series}

In this subsection we consider the action of a Hausdorff operator on ordinary Dirichlet series 
$$
D=\sum_{n=1}^\infty \frac{a(n)}{n^{s}}.
$$

Let $\Bbb{Z}^\infty$ be the additive group of   infinite
sequences of  integers with finite support, $\Bbb{Z}^\infty_+:=\{\alpha\in \Bbb{Z}^\infty: \forall k \alpha_k\ge 0\}$. Since every natural $n$ has the prime number decomposition $n=\mathfrak{p}^\alpha:=p_1^{\alpha_1}\dots p_N^{\alpha_N}$
where $\alpha\in \Bbb{Z}^\infty_+$, and $\mathfrak{p}=\{2,3,5,\dots\}$ the set of all primes, one can
identify the series $D$   with the corresponding coefficient function $\Bbb{Z}^\infty\backepsilon\alpha\mapsto a({\mathfrak{p}^\alpha})$
\textit{supported in} $\Bbb{Z}^\infty_+$.
 In this case the action of a Hausdorff operator on $D$ means the action  on the function $a({\mathfrak{p}^{(\cdot)}})$ and Definition 4 takes the form
 $$
 (\mathcal{H}_{\Phi,\tau}a({\mathfrak{p}^{(\cdot)}}))(\alpha)=\int_{\Omega}\Phi(u)a({\mathfrak{p}^{\tau_u(\alpha)}})d\mu(u).
 $$
(Since the function $a(\mathfrak{p}^{(\cdot)})$
is supported in $\Bbb{Z}^\infty_+$, one can consider only such automorphisms $\tau_u$ that $\tau_u(\alpha)\in \Bbb{Z}^\infty_+$.)

We show that a certain class of such operators
 acts in the Banach space $\mathcal{D}_\infty$ of all sums of ordinary Dirichlet series $D$
which converge and define a bounded and  holomorphic function  $D(\cdot)$ on the half-plane $\{\mathrm{Re} s > 0\}$ ($\mathcal{D}_\infty$  is endowed with the supremum norm $\|\cdot\|_\infty $  on $\{\mathrm{Re} s > 0\}$). We identify the function $D(\cdot)\in \mathcal{D}_\infty$   with the coefficient function $a({\mathfrak{p}^{(\cdot)}})$ as mentioned above and put $\|a({\mathfrak{p}^{(\cdot)}})\|:=\|D(\cdot)\|_\infty$.

\textbf{Theorem 6.} \textit{Let  $\Phi\in L^1(\mu)$, and a family $(\tau_u)_{u\in \Omega}$ of automorphisms of $\mathbb{Z}^\infty$ enjoys the property $\tau_u:(\Bbb{Z}^\infty_+)^c\to (\Bbb{Z}^\infty_+)^c$ a.e. $u\in \Omega$. Then a Hausdorff operator
$\mathcal{H}_{\Phi,\tau}$
acts in $\mathcal{D}_\infty$ and $\|\mathcal{H}_{\Phi,\tau}\|_{\mathcal{L}(\mathcal{D}_\infty)}\le \|\Phi\|_{L^1}$.
}

Proof.  The group $\Bbb{Z}^\infty$ can be identified with the dual of the infinite-dimensional torus $\mathbb{T}^\infty$  via the map $\alpha\mapsto \chi_\alpha$, where the character $\chi_\alpha(t)=t^\alpha:=t_1^{\alpha_1}\dots t_N^{\alpha_N}$ and $\alpha=(\alpha_1,\dots,\alpha_N,0,0,\dots)\in \Bbb{Z}^\infty$.

It is
proven in \cite{HLS} (see also \cite[Corollary 5.3]{DGMS} or  \cite[Theorem 6.2.3, p. 145]{QQ}) that  the map $\Psi$ that  takes a function $a({\mathfrak{p}^{(\cdot)}})$ from $\mathcal{D}_\infty$ to a function  $f_a$ on $\mathbb{T}^\infty$ with the Fourier transform $\widehat{f_a}(\alpha)=a({\mathfrak{p}^{\alpha}})$ ($\alpha\in \Bbb{Z}^\infty_+$)
 is an isometric isomorphism of Banach spaces   $\mathcal{D}_\infty$ and $H^\infty_{\Bbb{Z}^\infty_+}(\mathbb{T}^\infty)$.

Now Theorem 1 with  $G=\mathbb{T}^\infty$  shows that for $\alpha\in \Bbb{Z}^\infty_+$ one has
$$
(\mathcal{H}_{\Phi,(\tau^*)^{-1}}f_a)^{\wedge}(\alpha)=
(\mathcal{H}_{\Phi,\tau}\widehat{f_a})(\alpha)=(\mathcal{H}_{\Phi,\tau}a({\mathfrak{p}^{(\cdot)}}))(\alpha).
$$

Putting $A(u)=(\tau_u^*)^{-1}$, $p=\infty$ in Theorem 2 we get
$$
\mathcal{H}_{\Phi,(\tau^*)^{-1}}f_a=f_b,
$$
where $f_b\in  H^\infty_{\Bbb{Z}^\infty_+}(\mathbb{T}^\infty)$ and therefore $(\mathcal{H}_{\Phi,(\tau^*)^{-1}}f_a)^{\wedge}=\widehat{f_b}$. Since $\widehat{f_b}(\alpha)=b({\mathfrak{p}^{\alpha}})$ for all $\alpha\in \Bbb{Z}^\infty_+$, it follows that
$$
\mathcal{H}_{\Phi,\tau}a({\mathfrak{p}^{(\cdot)}})=b({\mathfrak{p}^{(\cdot)}}),
$$
i.e. $\mathcal{H}_{\Phi,\tau}$
acts in $\mathcal{D}_\infty$.

Finally,  for the isometric isomorphism  $\Psi:\mathcal{D}_\infty\to H^\infty_{\Bbb{Z}^\infty_+}(\mathbb{T}^\infty)$ we have $\Psi^{-1}f_a=a({\mathfrak{p}^{(\cdot)}})$ for each $f_a\in H^\infty_{\Bbb{Z}^\infty_+}(\mathbb{T}^\infty)$. So,
$$
\Psi\mathcal{H}_{\Phi,\tau}\Psi^{-1}f_a=\Psi\mathcal{H}_{\Phi,\tau}a({\mathfrak{p}^{(\cdot)}})=\Psi b({\mathfrak{p}^{(\cdot)}})=f_b.
$$
Thus $\Psi\mathcal{H}_{\Phi,\tau}\Psi^{-1}=\mathcal{H}_{\Phi,(\tau^*)^{-1}}$ and therefore

  $$\|\mathcal{H}_{\Phi,\tau}\|_{\mathcal{L}(\mathcal{D}_\infty)}=
 \|\mathcal{H}_{\Phi,(\tau^*)^{-1}}\|_{\mathcal{L}(H_{\Bbb{Z}^\infty_+}^\infty)}\le \|\Phi\|_{L^1}
 $$
  which completes the proof.

The next corollary is a generalization of a Theorem  of Bohr (see, e.g., \cite[p. 224]{Rud}).

\textbf{Corollary 5.} \textit{Let  $\Phi\in L^1(\mu)$, and a family $(\tau_u)_{u\in \Omega}$ of automorphisms of $\mathbb{Z}^\infty$ enjoys the property $\tau_u:(\Bbb{Z}^\infty_+)^c\to (\Bbb{Z}^\infty_+)^c$ a.e. $u\in \Omega$. Let $E$ be the set of  all $\alpha\in \Bbb{Z}^\infty_+$  with $\sum \alpha_j=1$.  Then for every $D(\cdot)\in \mathcal{D}_\infty$   with the coefficient function $a({\mathfrak{p}^{(\cdot)}})$ we have }

$$
\sum_{\alpha\in E}|(\mathcal{H}_{\Phi,\tau}a({\mathfrak{p}^{(\cdot)}}))(\alpha)|\le \|\Phi\|_{L^1}\|a({\mathfrak{p}^{(\cdot)}})\|.
$$

Proof. By Theorem 6 the function $\phi$ on $\{\mathrm{Re} s>0\}$ which is a sum of a  Dirichlet series with the coefficient function
$$
c({\mathfrak{p}^{(\cdot)}}):=\mathcal{H}_{\Phi,\tau}a({\mathfrak{p}^{(\cdot)}})
$$
belongs to $\mathcal{D}_\infty$. Then by Theorem  of Bohr mentioned above
$$
\sum_{\alpha\in E}|(\mathcal{H}_{\Phi,\tau}a({\mathfrak{p}^{(\cdot)}}))(\alpha)|=\sum_{\alpha\in E}|c({\mathfrak{p}^{\alpha}})|=\sum_{p\in \mathfrak{p}}|c(p)|\le\|\phi\|_{\infty}=\|c({\mathfrak{p}^{(\cdot)}})\|
$$
$$
=\|\mathcal{H}_{\Phi,\tau}a({\mathfrak{p}^{(\cdot)}})\|_{\infty}\le
\|\mathcal{H}_{\Phi,\tau}\|\|a({\mathfrak{p}^{(\cdot)}})\|\le \|\Phi\|_{L^1}\|a({\mathfrak{p}^{(\cdot)}})\|
$$
what was required.

Now we  consider the concrete family of automorphisms of $\mathbb{Z}^\infty$ that meet  the condition of Theorem 6.

Let $\Omega=\Bbb{Z}^\infty_+$ (with the counting measure). For each $u\in \Bbb{Z}^\infty_+$  define the map $\sigma_u:\Bbb{Z}^\infty\to \Bbb{Z}^\infty$ as follows
$$
\sigma_u(\alpha)=(\alpha_1, -u_1\alpha_1+\alpha_2,\dots, -u_{k-1}\alpha_{k-1}+\alpha_k,\dots).
$$
Then $\sigma_u\in\mathrm{ Aut}(\Bbb{Z}^\infty)$, and its inverse is given by the rule
$\sigma_u^{-1}(\beta):=\alpha$ where $\beta\in \Bbb{Z}^\infty$ and $\alpha\in \Bbb{Z}^\infty$ satisfies the following recurrent relation: $\alpha_1:=\beta_1$,  $\alpha_k:=\beta_k+u_{k-1} \alpha_{k-1}$ ($k\ge 2$).

\textbf{Corollary 6.} \textit{Let $\Phi\in\ell^1(\Bbb{Z}^\infty_+)$. Then a discrete  Hausdorff operator
$$
 (\mathcal{H}_{\Phi,\sigma}a({\mathfrak{p}^{(\cdot)}}))(\alpha)=\sum_{u, \sigma_u(\alpha)\in \Bbb{Z}^\infty_+}\Phi(u)a({\mathfrak{p}^{\sigma_u(\alpha)}})\eqno(4)
$$
acts in $\mathcal{D}_\infty$ and $\|\mathcal{H}_{\Phi,\sigma}\|_{\mathcal{L}(\mathcal{D}_\infty)}\le \|\Phi\|_{\ell^1}$.
}

Proof.
The automorphism $\sigma_u$ of $\mathbb{Z}^\infty$ do maps the set $(\Bbb{Z}^\infty_+)^c=\{\alpha\in \Bbb{Z}^\infty:  \exists k \alpha_k< 0\}$ into itself (indeed, if $\alpha_k$ is the first negative entry of $\alpha\in (\Bbb{Z}^\infty_+)^c$ and $\beta=\sigma_u(\alpha)$ then $\beta_k<0$). It remains to note that in  our case  the operator $\mathcal{H}_{\Phi,\sigma}$ has the form (4) because the function $a(\mathfrak{p}^{(\cdot)})$
is supported in $\Bbb{Z}^\infty_+$.

For another  result in this direction see  Corollary 7 below.

\

\subsection{General Dirichlet Series}

To formulate and prove similar results on general Dirichlet series we need some notation, definitions, and results from \cite{DS2019}, \cite{I. Schoolmann}, and \cite{DSservey}.

%show that some sort of Hausdorff operators act on the Banach space $\mathcal{D}_\infty$ of all ordinary Dirichlet series $D_a(s)=\sum_{n=1}^\infty a(n)/n^{s}$
%which converge and define a bounded and  holomorphic function on the half-plane $\{\mathrm{Re} s > 0\}$ (endowed with the supremum norm $\|\cdot\|_\infty $  on $\{\mathrm{Re} s > 0\}$).

Let $\lambda=(\lambda_n)$ be a non-negative strictly increasing sequence of real numbers tending to $\infty$ ("a frequency"). The value  $L(\lambda):=\limsup_{n\to\infty}(\log n)/\lambda_n$ (the maximal width of the
strip of convergence and non absolutely convergence of the corresponding Dirichlet series) is  associated to a frequency $\lambda$.

A compact Abelian group $G$ is called \textit{a  $\lambda$-Dirichlet group} if there is a continuous
homomorphism $\beta:\mathbb{ R}\to G$ with dense range such that every continuous  character $\widehat{\lambda_n}=e^{-i \lambda_n\cdot}$ of $\mathbb{R}$ has an "extension" $h_{\lambda_n}\in X$  (which then is unique) such that  $h_{\lambda_n}\circ \beta=\widehat{\lambda_n}$.

We consider formal general Dirichlet series
$$
D_{\lambda}=\sum_{n=1}^\infty a_n(D)e^{-\lambda_n s}.
$$

In \cite{DS2019} the next two spaces were  introduced
$$
\mathcal{D}_\infty(\lambda):=\{D_{\lambda}: D_{\lambda} \mbox{ converge to a
  function from } H^\infty(\{\mathrm{Re} > 0\})\},
$$
and 
$\mathcal{D}_\infty^{ext}(\lambda)$
of all somewhere convergent $\lambda$-Dirichlet series, which have a holomorphic
and bounded extension to the right half-plane $\{\mathrm{Re} > 0\}$. In general $\mathcal{D}_\infty(\lambda)\subseteq \mathcal{D}_\infty^{ext}(\lambda)$ and Theorem 2.2 from \cite{DS2019} gives sufficient conditions for the equality here. Moreover, if $L(\lambda)<\infty$
then the space $\mathcal{D}_\infty^{ext}(\lambda)$ is complete with respect to the supremum norm over $\{\mathrm{Re} > 0\}$ \cite[Theorem 5.1]{I. Schoolmann}.

Let $(G,\beta)$ be a  $\lambda$-Dirichlet group. Following \cite{DS2019} for $f\in L^1(G)$ we consider formal general Dirichlet series of the form
$$
D_{f,\lambda}=\sum_{n=1}^\infty \widehat{f}(h_{\lambda_n})e^{-\lambda_n s}.\eqno(5)
$$
If the space $\mathcal{D}_\infty^{ext}(\lambda)$ is complete one has 
$$
\mathcal{D}_\infty^{ext}(\lambda)=\mathcal{D}_\infty(\lambda)=\{D_{f,\lambda}: f\in H^\infty_E(G) \mbox{ where  } E=\{h_{\lambda_n}: n\in \mathbb{N}\}\}\eqno(6)
$$
(see \cite[Theorem 4.1]{DSservey} and references therein).

We introduce a \textit{Hausdorff operator on (formal) general Dirichlet series } of the form (5) as follows.

\textbf{Definition 6.} Let $(G,\beta)$  be a  $\lambda$-Dirichlet group,  $\Phi\in L^1(\mu)$, and $\{\tau_u:u\in \Omega\}\subset \mathrm{Aut}(X)$. For $f\in L^1(G)$ we put
$$
\mathrm{H}_{\Phi,\tau}D_{f,\lambda}:=D_{g,\lambda},
$$
where $g=\mathcal{H}_{\Phi,(\tau^*)^{-1}}f$.

(This definition is correct, because $g\in L^1(G)$.)

Since $\widehat{g}= \mathcal{H}_{\Phi,\tau} \widehat{f}$ by  Theorem 1,
Definition 6 means that 
$$
\mathrm{H}_{\Phi,\tau}:\sum_{n=1}^\infty \widehat{f}(h_{\lambda_n})e^{-\lambda_n s}\mapsto \sum_{n=1}^\infty (\mathcal{H}_{\Phi,\tau} \widehat{f})(h_{\lambda_n})e^{-\lambda_n s}.
$$

\textbf{Theorem 7.}\textit{ Let  $(G,\beta)$  be a  $\lambda$-Dirichlet group, $E:=\{h_{\lambda_n}:n\in \mathbb{N}\}$, $\tau_u:E^c\to E^c$, and $\Phi\in L^1(\mu)$. If $L(\lambda)<\infty$ then
$\mathrm{H}_{\Phi,\tau}$ acts in $\mathcal{D}_\infty(\lambda)$ and $\|\mathrm{H}_{\Phi,\tau}\|_{\mathcal{L}(\mathcal{D}_\infty)}\le \|\Phi\|_{L^1}$.
}

Proof. Since $L(\lambda)<\infty$, $\mathcal{D}_\infty^{ext}(\lambda)$ is complete by \cite[Theorem 5.1]{I. Schoolmann} and therefore (6) holds.   Let $D_{f,\lambda}\in\mathcal{D}_\infty(\lambda)$. Then $f\in H^\infty_E(G)$ and the function $g=\mathcal{H}_{\Phi,(\tau^*)^{-1}}f$ belongs to the space $H^\infty_E(G)$, too,
by Theorem 2. Thus, the operator $\mathrm{H}_{\Phi,\tau}$ acts in $\mathcal{D}_\infty(\lambda)$.

Following \cite{DS2019} consider the Bohr map
$$
\mathcal{B}:H^\infty_E(G)\to \mathcal{D}_\infty(\lambda),\ f\mapsto D_{f,\lambda}.
$$
 As was mentioned above, in our case $\mathcal{D}_\infty^{ext}(\lambda)=\mathcal{D}_\infty(\lambda)$. So, Theorem 4.12 from \cite{DS2019} states that $\mathcal{B}$ is an  isometrical isomorphism of this Banach spaces. But the equality $\mathrm{H}_{\Phi,\tau}D_{f,\lambda}=D_{g,\lambda}$ means that $\mathrm{H}_{\Phi,\tau}\mathcal{B}f=\mathcal{B}g=\mathcal{B}\mathcal{H}_{\Phi,(\tau^*)^{-1}}f$ for all $f\in H^\infty_E(G)$. In other words,
$$
\mathrm{H}_{\Phi,\tau}=\mathcal{B}\mathcal{H}_{\Phi,(\tau^*)^{-1}}\mathcal{B}^{-1}.
$$
It follows that $\|\mathrm{H}_{\Phi,\tau}\|_{\mathcal{L}(\mathcal{D}_\infty)}=\|\mathcal{H}_{\Phi,(\tau^*)^{-1}}\|_{\mathcal{L}(H^\infty_E)}\le \|\Phi\|_{L^1}$. This completes the proof.

One can apply Theorem  7 to the space $\mathcal{D}_\infty=\mathcal{D}_\infty((\log n))$ of ordinary Dirichlet series and get the next

\textbf{Corollary 7.} \textit{Let  $G=\mathbf{b}\mathbb{R}$ be the Bohr compactum, $\Omega=\{1/q:q\in \mathbb{N}\}$, and $\tau_u(\gamma)=u\gamma$ for $u\in \Omega$, $\gamma\in \mathbb{R}$. If $\Phi\in \ell^1(\Omega)$ then the discrete Hausdorff operator $\mathrm{H}_{\Phi,\tau}$ acts in the space  $\mathcal{D}_\infty$ of ordinary Dirichlet series and $\|\mathrm{H}_{\Phi,\tau}\|_{\mathcal{L}(\mathcal{D}_\infty)}\le \|\Phi\|_{L^1}$.}

Proof. First note that by \cite[Example 3.19]{DS2019} the Bohr compactum $(\mathbf{b}\mathbb{R},\beta)$
where $\beta(t)=\widehat{t}$ is a   $\lambda$-Dirichlet group for any frequency $\lambda$. Further,  in our case $\lambda_n=\log n$. If we identify the group $\mathbb{R}_d$ with the dual for $\mathbf{b}\mathbb{R}$ then   $E=\{h_{\lambda_n}:n\in \mathbb{N}\}=\{\log n: n\in \mathbb{N}\}$. 
 Since $\tau_u: E^c\to E^c$ for all $u\in\Omega$ and
 $L((\log n))=1$, the result follows from Theorem 7.

\

\section{Examples  in the Case of Ordered Dual}

\textbf{Example 3.} Let $G=\mathbf{b}\mathbb{R}$ be the Bohr compactification of the reals, $X=\mathbb{R}_d$  as in Example 2.
 The map $\tau_u:X\to X, \gamma\mapsto u\gamma$ belongs to $\mathrm{Aut}_+(X)$ for $u\in (0,\infty)$. Since $(\tau_u^*)^*=\tau_u$ \cite[(24,41)]{HiR}, it follows that the map $A(u):=\tau_u^*$  belongs to $\mathrm{Aut}(\mathbf{b}\mathbb{R})^+$ for each $u>0$. If $\mu$ is some regular Borel measure on $\Omega:=(0,\infty)$, the corresponding Hausdorff operator on $\mathbf{b}\mathbb{R}$ is
$$
 \mathcal{H}_{\Phi,\tau^*}g(x)=\int_{(0,\infty)} \Phi(u)g(\tau_u^*(x))d\mu(u),\ \ x\in \mathbf{b}\mathbb{R}
$$
(recall that $\tau_u^*(x)=x\circ \tau_u$). This operator is bounded on $H^p(\mathbf{b}\mathbb{R})$ ($1\le p\le\infty$),   $BMOA(\mathbf{b}\mathbb{R})$, $H^1_\mathbb{R}(\mathbf{b}\mathbb{R})$ (for real valued $\Phi$), and   $BMO(\mathbf{b}\mathbb{R})$  if and only if $\Phi\in L^1(\mu)$ and its norm does not exceed $\|\Phi\|_{L^1(\mu)}$.

\textbf{Example 4.} Let $G=\mathbb{T}^d$ be the  $d$-dimensional torus ($d\ge 2$).  Let $\Omega$ be the subgroup of the arithmetic  group $\mathrm{GL}(d,\mathbb{Z})$ which consists of  matrices $u=(u_{ij})_{i,j=1}^d$ with $\det u=\pm 1$. Then every map
 $$
A(u)(z)= z^u:=(z_1^{u_{11}}z_2^{u_{12}}\dots z_d^{u_{1d}},\dots,z_1^{u_{d1}}z_2^{u_{d2}}\dots z_d^{u_{dd}})
$$
($z=(z_j)_{j=1}^d\in \mathbb{T}^d$) belongs to $\mathrm{Aut}(\mathbb{T}^d)$  (see, e.g., \cite[(26.18)(h)]{HiR}).
Thus, the corresponding Hausdorff operator over $\mathbb{T}^d$ takes the form
$$
(\mathcal{H}_{\Phi,A}f)(z)=\int_{\Omega}\Phi(u)f(z^u)d\mu(u)
$$
where $\mu$ stands for some regular Borel measure on $\Omega$ (e.~g., $\mu$ is a  Haar measure of the  group  $\Omega$).

Every character of $ \mathbb{T}^d$ has the form $\chi_n(z)=z_1^{n_1}\dots z_d^{n_d}$, where
$n=(n_1,\dots,n_d)\in\mathbb{Z}^d$. Thus, the dual of $ \mathbb{T}^d$ can be identified with the group $\mathbb{Z}^d$ via the map $\chi_n\mapsto n$. We endow $\mathbb{Z}^d$ with the lexicographic order. For this order the positive cone is
  $$
 X_+=\{n\in \mathbb{Z}^d: n_1>0\}\cup\{n\in \mathbb{Z}^d: n_1=0, n_2>0\}\cup
 $$
$$
  \cdots \cup \{n\in \mathbb{Z}^d: n_1= n_2=\ldots=n_{d-1}=0,   n_d>0\}\cup\{0\}.
 $$

 Consider the arithmetic strict lower triangular group $\mathrm{T}_1(d,\mathbb{Z})$. This group consists of   matrices $u\in\mathrm{SL}(d,\mathbb{Z})$ such that $u_{ii}=1$, and $u_{ij}=0$ for $i<j$. Then the map
 $$\tau_u(n):= un^\top=(n_1, u_{21}n_1+n_2,\dots, u_{d1}n_1+\dots + u_{d,d-1}n_{d-1}+n_d)
 $$
  (here $n^\top\in\mathbb{Z}^d$ is a column vector) belongs to  $\mathrm{Aut}_+(X)$. Since
 $$
 \tau_u^*(z)=(z_1,z_1^{u_{21}}z_2,\dots, z_1^{u_{d1}}z_2^{u_{d2}}\dots z_{d-1}^{u_{d,d-1}}z_d),
  $$
  in this case,
 $$
 (\mathcal{H}_{\Phi,\tau^*}f)(z)=\int_{\mathrm{T}_1(d,\mathbb{Z})}\Phi(u)f(z_1,z_1^{u_{21}}z_2,\dots, z_1^{u_{d1}}z_2^{u_{d2}}\dots z_{d-1}^{u_{d,d-1}}z_d)d\mu(u)
 $$
where $\mu$ is some regular Borel measure on $\mathrm{T}_1(d,\mathbb{Z})$.

This operator is bounded on $H^p(\mathbb{T}^d)$ ($1\le p\le\infty$),   $BMOA(\mathbb{T}^d)$, $H^1_\mathbb{R}(\mathbb{T}^d)$ (for real valued $\Phi$),  and  $BMO(\mathbb{T}^d)$  if and only if $\Phi\in L^1(\mu)$ and its norm does not exceed $\|\Phi\|_{L^1(\mu)}$.

\textbf{Example 5.} Let $\mathbb{T}^\infty$ be the infinite-dimensional torus and $X=\Bbb{Z}^\infty_{{\rm lex}}$ --- the additive group of   infinite
sequences of  integers with finite support endowed with the  lexicographic order. For this order, by  definition  the positive cone is
$$
   X_+=\{0\}\cup\{\alpha\in \mathbb{Z}^\infty: \alpha_1>0\}\cup\{\alpha\in \mathbb{Z}^\infty: \alpha_1=0, \alpha_2>0\}
   \cup\dots.
$$
In other words, $X_+$ consists of sequences whose first non-zero entry is positive and the  zero sequence.

 As above we  identify the group $\Bbb{Z}^\infty$   with the dual group of  $\mathbb{T}^\infty$  via the map $\alpha\mapsto \chi_\alpha$ where  $\chi_\alpha(z)=z_1^{\alpha_1}z_2^{\alpha_2}\dots$ ($z\in \mathbb{T}^\infty$).

Let $\mathrm{J}(\infty,\mathbb{Z})$ consists of  infinite lower two-diagonal matrices $u$ of integers such that $u_{ii}=1$, $u_{ij}=0$ for $i<j$, and $u_{k,1}=\dots=u_{k,k-2}=0$ for $k\ge 3$.

Then the map
$$
\tau_u(\alpha):= u\alpha^\top= (\alpha_1, u_{21}\alpha_1+\alpha_2,u_{32}\alpha_2+\alpha_3,\dots, u_{k,k-1}\alpha_{k-1}+\alpha_k,\dots)
$$
($u\in \mathrm{J}(\infty,\mathbb{Z})$,  $\alpha\in\mathbb{Z}^\infty$) belongs to  $\mathrm{Aut}_+(\Bbb{Z}^\infty_{{\rm lex}})$. Since
 $$
 \tau_u^*(z)=z^u:=(z_1,z_1^{u_{21}}z_2,\dots, z_{k-1}^{u_{k,k-1}}z_k,\dots),
  $$
  in this case,
 $$
 (\mathcal{H}_{\Phi,\tau^*}f)(z)=\int_{\mathrm{J}(\infty,\mathbb{Z})}\Phi(u)f(z^u)d\mu(u)
 $$
where $\mu$ is some regular Borel measure on $\mathrm{J}(\infty,\mathbb{Z})$.

This operator is bounded on $H^p(\mathbb{T}^\infty)$ ($1\le p\le\infty$),   $BMOA(\mathbb{T}^\infty)$, $H^1_\mathbb{R}(\mathbb{T}^\infty)$ (for real valued $\Phi$),  and  $BMO(\mathbb{T}^\infty)$  if and only if $\Phi\in L^1(\mu)$ and its norm does not exceed $\|\Phi\|_{L^1(\mu)}$.

\textbf{Example 6.} Let $\mathbf{a}=(2,3,4,\dots)$. Then \textit{the $\mathbf{a}$-adic solenoid} $\Sigma_\mathbf{a}$  (see, e.g., \cite[(10.12)]{HiR}) is a compact and connected topological group and is topologically isomorphic to  the  character
group  $\widehat{\mathbb{Q}_d}$ of the discrete additive group $X=\mathbb{Q}_d$ of rationals
\cite[(25.4)]{HiR}. On the other hand, by \cite[(25.5)]{HiR} the group $\widehat{\mathbb{Q}_d}$
can be identified with some  subgroup $G$ of  the infinite-dimensional torus  $\mathbb{T}^\infty$ in the following way. Let the sequence  $\alpha=(\alpha_n)_{n\in\mathbb{ N}}\in \mathbb{T}^\infty$ be such that $\alpha_n=\alpha_{n+1}^{n+1}$ for all $n\in\mathbb{ N}$. Then it produces a character
of   $\mathbb{Q}_d$ via the rule
$$
\chi_\alpha\left(\frac{m}{n!}\right)=\alpha_{n}^{m} \ \ (m\in \mathbb{Z}, n\in \mathbb{N}).
$$
Moreover, each character of $\mathbb{Q}_d$ can be identified with such a sequence $\alpha$
and we get an isomorphism $\alpha\mapsto \chi_\alpha$ of the subgroup  $G:=\{\alpha\}\subset \mathbb{T}^\infty$ and $\Sigma_\mathbf{a}$.
Thus, one can identify the group $G$ with  $\Sigma_\mathbf{a}$.  Further, for each $q\in \mathbb{Q}$, $q>0$ the map $l_q(x)=qx$ is  an order automorphism of the group $\mathbb{Q}_d$ endowed with the usual order. It  follows that the corresponding dual automorphism $l_q^*$ of the dual group  $G=\Sigma_\mathbf{a}$ belongs to $\mathrm{Aut}(\Sigma_\mathbf{a})^+$. This yields  that for every measurable map $k:\Omega\to \mathbb{Q}_+\setminus\{0\}$   the corresponding  Hausdorff operator
$$
\mathcal{H}_{\Phi,l_k^*}f(\alpha)=\int_\Omega \Phi(u)f(l_{k(u)}^*(\alpha))d\mu(u)
$$
is bounded on  $H^p(\Sigma_\mathbf{a})$ ($1\le p\le\infty$),   $BMOA(\Sigma_\mathbf{a})$, $H^1_\mathbb{R}(\Sigma_\mathbf{a})$ (for real valued $\Phi$),  and  $BMO(\Sigma_\mathbf{a})$  if and only if $\Phi\in L^1(\mu)$ and its norm does not exceed $\|\Phi\|_{L^1(\mu)}$.

\textbf{Example 7.} Let $G$ be  a  compact and connected Abelian group with totally ordered dual and $\Omega$  a compact subgroup of $\mathrm{Aut}(G)$ with normalized Haar measure $\mu$.
 The \textit{generalized   shift operator of Delsarte}  \cite{Delsarte}, \cite[Ch. I, \S 2]{Lev} (also  the  terms “generalized  translation operator of Delsarte”, or “generalized displacement operator of Delsarte” are used) is defined to be
$$
T^hf(x)=\int_{\Omega} f(hu(x))d\mu(u)\quad (x,h\in G).
$$
Then $T^{h}=\mathcal{H}_1S_h$, where
$$
\mathcal{H}_1f(x):=\int_\Omega f(u(x))d\mu(u)
$$
is a Hausdorff operator on $G$ with $\Phi\equiv 1$,  $A(u)=u$, and $S_h f(x):=f(hx)$.
   Let $u\in \mathrm{Aut}(G)^+$ for $\mu$-a.~e. $u\in \Omega$. Then for every fixed $h$  the generalized   shift operator of Delsarte is bounded on $H^p(G)$ ($1\le p\le\infty$), $BMOA(G)$, $BMO(G)$, and $H^1_{\mathbb{R}}(G)$. In addition, its norm in this spaces  equals  to $\mu(\Omega)=1$    (Remark 1).

\

\textbf{Acknowledgments}.   This work was supported by the State Program of Scientific Research of the Republic of Belarus.

\

\end{document}